\documentclass{amsart}

\newtheorem{thm}{Theorem}[section]

\newtheorem{DFN}[thm]{Definition}
\newtheorem{RMK}[thm]{Remark}
\newtheorem{EG}[thm]{Example}

\newcommand{\cal}{\mathcal}
\newcommand{\C}{{\mathbb C}}
\newcommand{\PP}{{\mathbb P}}
\newcommand{\Z}{{\mathbb Z}}
\newcommand{\pf}{\par \noindent {\em Proof.}\hspace{.5em}}
\newcommand{\edpf}{ \framebox(5,7)[bl]{} \par \bigskip \par}
\newcommand{\Sing}{\mbox{\rm Sing\,}}
\newcommand{\St}{\mbox{\it St}}
\newcommand{\Pic}{\mbox{\rm Pic}\;}
\newcommand{\im}{\mbox{\rm im}\;}
\newcommand{\ie}{\textit{i.e.\ }}
\newcommand{\cf}{\textit{cf.\ }}
\newcommand{\yuan}{\circle{2}}

\begin{document}

\title{The Embedded Resolution of $f(x,y)+z^2: (\C^3,0)\to (\C,0)$}
\author{Chunsheng Ban}
\email{cban@math.ohio-state.edu}
\author{Lee J. McEwan}
\email{mcewan@math.ohio-state.edu}
\author{Andr\'as\ N\'{e}methi}
\email{nemethi@math.ohio-state.edu}
\address{Department of Mathematics, Ohio State University}

\maketitle

\section{Introduction}

\subsection{}
The goal of the present paper is the presentation of an
``embedded resolution'' of $\{f(x,y)+z^2=0\},0)\subset (\C^3,0)$
using the method of Jung. In the first part of the introduction, we present
the terminology and the strategy of the paper.

Let $(Y,0)$ be the germ of an analytic space. We say that
a proper, surjective analytic map $\phi: \widetilde{Y} \to U$ is a
resolution of $(Y,0)$ (where $U$ is a small representative of $(Y,0)$),
if $\widetilde{Y}$ is smooth,
$\phi^{-1}(U-\Sing Y)$ is dense in $\widetilde{Y}$,
 and $\phi^{-1}(U-\Sing Y)\to U-\Sing Y$ is an
isomorphism.

More generally, if $(Y,0)$ is a local divisor in $(\C^n,0)$, we say
that $\phi: \widetilde{X}\to U$ is an embedded resolution of the
pair $(Y,0)\subset (\C^n,0)$ (where again $U$ is a small representative
of $(\C^n,0)$), if $\tilde{X}$ is smooth, $\phi$ is an isomorphism
above $U-\Sing Y$, and $\phi^{-1}(Y)$ is a normal crossing divisor
in $\widetilde{X}$.

Jung's strategy (\cite{Jung}, see also \cite{Lipman1}, \cite{Lipman2})
gives a recipe how
one can reduce the construction of the resolution of the $n$-dimensional
space $(Y,0)$ to the case of the construction of an embedded resolution
of (the ``lower dimensional case'') $(\Delta,0)\subset (\C^n,0)$ and
the resolution of the so-called quasi-ordinary singularities of
dimension $n$. (For quasi-ordinary singularities, see \cite{Lipman1},
\cite{Lipman2}). Indeed, consider a projection
$p:(Y,0)\to (\C^n,0)$ with finite fibers,
whose reduced discriminant locus is $(\Delta,0)$. Then resolve the pair
$(\Delta,0)\subset (\C^n,0)$ by $\phi:Z\to U\subset\C^n$
and pull back $p$ to $p':Y'\to Z$. Then
$Z$ is smooth and the branch locus of $p'$ is a normal
crossing divisor, hence $Y'$ has only quasi-ordinary
singularities. Resolving these singularities  we obtain $\widetilde{Y}$.

For example, consider the hypersurface singularity
$(Y,0)=(\{f(x,y)+z^2=0\},0)$, where $f$ is an isolated plane curve
singularity. Then $p:(Y,0)\to (\C^2,0)$, given by
$(x,y,z)\mapsto (x,y)$ is a double covering with branch locus
$(\{f=0\},0)\subset (\C^2,0)$. Using the above notations, if
$\phi: Z\to U\subset\C^2$ is an embedded resolution of
this plane curve singularity, then $Y'$ normalized has
only Hirzebruch-Jung type singularities of type $\{z^2=x^my^n\}^{norm}$.

If  we want to obtain an \textit{embedded} resolution by  a similar strategy,
then we face two problems. To explain this, start with
$(Y,0)\subset(\C^n,0)$ and a projection
$p:(\C^n,0)\to (\C^{n-1},0)$ such that $p|Y $ is finite.
Let $(\Delta,0)$ be the discriminant
locus of $p|Y$. Then take an embedded resolution
$\phi:Z\to U^{n-1}\subset\C^{n-1}$ of the pair
$(\Delta,0)\subset (\C^{n-1},0)$ and construct by pull-back the diagram

\[
\begin{array}{rcl}
X & \stackrel{\phi'}{\longrightarrow} &   U^n\subset \C^n \\
\vcenter{\llap{\scriptsize $p'$}}\Big\downarrow & &
  \Big\downarrow\vcenter{\rlap{\scriptsize $p$}} \\
Z & \stackrel{\phi}{\longrightarrow} &   U^{n-1} \subset \C^{n-1}
\end{array}
\]

\noindent Then still we have to find the embedded resolution of the
pair $(\phi')^{-1}(Y)\subset X$. Although $(\phi')^{-1}(Y)$
has only quasi-ordinary singularities, in general we know
very little about their {\em embedded} resolution.

Moreover,   even if we solve this problem (\ie we find an embedded resolution
$\widetilde{X}\to X$ of these quasi-ordinary singularities),
the final result of the above
construction has a small beauty defect: the constructed birational
map $\widetilde{X}\to U^n$ is not an isomorphism above the
complement of $\Sing Y$, but only above a smaller set, the
complement of $p^{-1}(\Sing\Delta)$.

For example, if $(Y,0)=(\{f(x,y)+z^2=0\},0)\subset(\C^3,0)$,
and $p(x,y,z)=(x,y)$, and we solve all the technical problems in the
above program, then we obtain a map
$\widetilde{\phi}:\widetilde{X}\to U^3\subset\C^3$ so that
$\widetilde{\phi}^{-1}(Y)$ is a normal crossing divisor, but
$\widetilde{\phi}$ is an isomorphism only above $\C^3-\{x=y=0\}$,
not above $\C^3-\{0\}$.

Even if this birational modification does not cover exactly the above
definition of the embedded resolution, in almost all the applications
it plays the role of an embedded  resolution:
all the invariants which can be read from an embedded resolution
can be read from this modification as well.

The goal of the present paper is exactly the presentation of {\em this}
birational modification (``embedded resolution'') of
$(\{f(x,y)+z^2=0\},0)\subset (\C^3,0)$.

\subsection{Preliminary Remarks}
Let $g:(\C^n,0)\to (\C,0)$ be the germ of an analytic function.
Sometimes,  an embedded resolution
$\widetilde{\phi}: \widetilde{X}\to U$
of the pair $(\{g=0\},0)\subset(\C^n,0)$
(or a birational modification as above)
is called the ``resolution of $g$''.
$E$ denotes the exceptional divisor, and  $E=E_1\cup \cdots \cup  E_s$
the decomposition  of $E$  into its irreducible components.

In our situation when $g=f(x,y)+z^2$,  since $\widetilde{\phi}$
is an isomorphism above $\C^3-\{x=y=0\}$,
the compact irreducible exceptional
divisors are situated above the origin (\ie $\widetilde{\phi}(E_i)=0$),
while the non-compact irreducible exceptional divisors projects via
$\widetilde{\phi}$ onto the disc $\{x=y=0\}$.

When we want to codify the modification $\widetilde{\phi}$, we have to
decide what kind of information we would like to extract from
it. In general $\widetilde{\phi}:\widetilde{X}\to U$ (or even $E$)
carries a lot of analytic information which is impossible to codify in
any topological, numerical
 or combinatorial object. But if we want to study the
pair $(g^{-1}(0),0)\subset (\C^n,0)$ from a topological point of view,
it is enough to record only the
topological/numerical/combinatorial invariants of
$\widetilde{\phi}$. If we
want to codify additional analytic information as well, then we face the
problem of analytic (or algebraic) classification of varieties and
vector bundles (which is basically an unsolved problem). Therefore, in
general, we have to find the right compromise which is still
satisfactory for our final goals.

Take, for example the case $n=2$ and an arbitrary embedded resolution
$\widetilde{\phi}$. All the topology is completely
codified in the dual resolution graph of $E\subset\widetilde{X}$ with
the decorations $\{E_i\cdot E_i\}_i$ (the self-intersections in
$\widetilde{X}$) and the vanishing orders $\{m_i\}_i$ of
$g\circ\widetilde{\phi}$ along $E_i$'s. Indeed, from this decorated
graph the homeomorphism type of $(\widetilde{X},E)$,  or of
$(g^{-1}(0),0)\subset (\C^n,0)$,  can be completely recovered by
plumbing. If we want to recover the analytic type of the reduced exceptional
divisor $E$, then the combinatorics of the graph cannot identify the
position of the points $\bigcup_{j\ne i}E_j\cap E_i$ on $E_i$.
Therefore, in this second level, if we wish to recover the isomorphic
type of $E$ from our codification, we need also to keep in this
codification the pairs $(E_i,\bigcup_{j\ne i}E_j\cap E_i)_i$.
On the other hand, at the
third level, if one wants to codify all the analytic information about
$(\widetilde{X},E)$, then that is a rather difficult  problem: the
underlying topology, in general, carries many analytic structures which
are parametrized by a rather complicated moduli space.

In higher dimensions, the problem is more complicated, and in general
it is not clear at all what the good, convenient levels and codifications
of the resolution are.

In the situation discussed in this article,
when $g(x,y,z)=f(x,y)+z^2: (\C^3,0)\to(\C,0)$, it turns out that the
codification of the modification of $g$ constructed above
is closely related to the
codification of the embedded resolution of $f$. We will define the
analog of the resolution graph of $f$ for $g$, which will keep all the
topological information about the pair ($E\subset\widetilde{X}$) up
to a homeomorphism. This graph will not identify $E$ modulo an analytic
isomorphism, but the ambiguity will be very similar to the plane curve
singularity case. (In other words, $g$ carries the same amount of
analytic information as $f$.)

Actually,  our modification $\widetilde{\phi}=\phi_g$
of $g$ will be constructed
from a fixed embedded resolution $\phi=\phi_f$ of $f$. If we denote the
irreducible exceptional divisors of $\phi_f$ by $\{A_i\}_i$, then any
compact irreducible exceptional divisor $E_k$ of $\phi_g$ will be a
(possibly non-minimal) ruled surface with natural projection
$E_k\to A_i=\PP^1$. All the special fibers will be situated over
the intersection points $\{A_j\cap A_i\}_{j\ne i}\subset A_i$.
Similarly, the non-compact irreducible exceptional divisors will be
non-minimal disc bundles over some curves  $A_i$. If from the resolution
of $\phi_f$ we retain the information of the position of the intersection
points $\{A_j\cap A_i\}_{j\ne i}\subset A_i$, then the analytic type of
the exceptional divisors
$E_k$ will be completely determined. If we use only the dual resolution
graph of $f$, then in the analytic type of $E_k$ we will have the
ambiguity of the position of the special fibers. This
ambiguity will disappear in the computation of
any kind of numerical invariant, and in the identification
of different elements (like
$K_{E_k}$ or  $N_{\widetilde{X}|E_k}$)  in $\Pic(E_k)$ for any compact $E_k$.

\section{Review of Ruled Surfaces over $\PP^1$}

\noindent For a general reference for ruled surfaces, we recommend
\cite{Hartshorne}.

\subsection{}
Any ruled surface over a smooth curve is obtained as $\PP(\cal{E})$ of a
locally free sheaf (vector bundle) $\cal{E}$ of rank 2. Actually,
$\PP(\cal{E})=\PP(\cal{E}\otimes\cal{L})$ for any line bundle $\cal{L}$.
But, over $\PP^1$, any $\cal{E}$ can be written as
$\cal{E}=\cal{O}(a)\oplus\cal{O}(b)$, hence any ruled surface
over $\PP^1$ has the
``normal form'' $X_e=\PP(\cal{O}\otimes\cal{O}(-e))$ for some integer
$e\ge 0$. The surface $X_e$ can be obtained by gluing two copies of
$\C\times\PP^1$ (with coordinates $(x,[u_0:u_1])$ and $(y,[w_0:w_1])$
respectively) along $\C^{\ast}\times\PP^1$ by the identifications
$y=1/x$, and $[w_0:w_1]=[u_0:x^eu_1]$. One has a natural projection
$\pi:X_e\to \PP^1$ (in coordinates $(x,[u_0:u_1])\mapsto x$,
$(y,[w_0:w_1])\mapsto y$). Here $\PP^1=\C\bigsqcup_{\C^{\ast}}\C$, where
the $x$-chart $\C$ corresponds to $\{[\alpha:\beta]\mid\alpha\ne 0\}$
with $x=\beta/\alpha$.

An automorphism $\phi:X_e\to X_e$ with $\pi\circ\phi=\pi$ is called
a $\pi$-automorphism; their collection  form the group $\cal{G}(X_e,\pi)$.

The projection $\pi:X_e\to\PP^1$ has two natural sections with images
$C_0$ and $C_1$. $C_0$ is given by $\{u_1=w_1=0\}$ and has self
intersection $C_0^2=-e$; $C_1$ is given by $\{u_0=w_0=0\}$ and has self
intersection $C_1^2=e$. Obviously $C_0\cap C_1=\emptyset$.

\subsection{Facts}

{\em
a) (\cite{Hartshorne}, V. 2.3)
$\Pic X_e=\mathbb{Z}\oplus\pi^{\ast}\Pic \PP^1=\mathbb{Z}\oplus\mathbb{Z}$
is generated by $C_0$ and an arbitrary fixed fiber $f$ of $\pi$. They
satisfy $C_0\cdot f=1$, $f^2=0$.

b) (\cite{Hartshorne}, V 2.8)
The  invertible sheaf $\cal{O}_{X_e}(1)$ (provided
by the Proj-construction of $\PP(\cal{E})$) satisfies
$\cal{O}_{X_e}(1)=\cal{O}_{X_e}(C_0)$.
}

\subsection{}
If $\im\tau$ is the image of a section $\tau:\PP^1\to X_e$, then
$\im\tau\equiv C_0+nf$ in $\Pic X_e$ for some $n$. The condition
$\im\tau\cap C_0=\emptyset$ implies $(C_0+nf)\cdot C_0=0$, hence $n=e$. In
particular $(\im\tau)^2=(C_0+ef)^2=e$.

\subsection{}
If $e=0$, then $X_0$ can be represented as
$\pi:X_0=\PP^1\times\PP^1\to\PP^1$,  where $\pi$ is the first
projection and $C_0$ and $C_1$ are two fibers of the projection on the
second factor. Obviously, the system $(\pi:X_0\to\PP^1;C_0,C_1)$ is
uniquely determined modulo the action of $\cal{G}(X_0,\pi)$
(i.e. it does not depend  on the choice of $C_0$ and $C_1$).

\subsection{}
Assume that $e>0$. In the sequel we will prove that the last sentence
of (2.4) is valid in this case as well. First notice that the section
$C_0$ is uniquely determined by the condition $C_0^2=-e$
(\cf \cite{Hartshorne}, V 2.11.3). Therefore, Fact 2.2 implies that any
$\phi\in\cal{G}(X_e,\pi)$ keeps $C_0$, $\cal{O}_{X_e}(1)$, $\Pic X_e$
invariant.

By \cite{Hartshorne}, V 2.6, there is a one-to-one correspondence between
sections $\tau:\PP^1\to X_e$ of $\pi$ and surjections
$\cal{O}\oplus\cal{O}(-e)\to\cal{L}$, where
$\cal{L}\in\Pic\PP^1$, given by $\cal{L}=\tau^{\ast}\cal{O}_{X_e}(1)$.
Under this correspondence, $C_0$ corresponds to the second projection
${\it pr}_2:\cal{O}\oplus\cal{O}(-e)\to\cal{O}(-e)$ and $C_1$ to
the first projection ${\it pr}_1:\cal{O}\oplus\cal{O}(-e)\to\cal{O}$.

Notice that it is possible to construct sections $\tau:\PP^1\to X_e$
with image $C_1':=\tau(\PP^1)$, not identical to $C_1$, but satisfying
the same numerical properties as $C_1$:  $C_1'\cap C_0=\emptyset$ and
${C_1'}^2=e$. Nevertheless, one has:

\subsection{Proposition} {\em
Fix $e>0$ and consider a ruled surface $\pi: X_e\to\PP^1$ and $C_0$,
$C_1\subset X_e$ as above. Then take an arbitrary section $\tau$ with
$\im\tau=C_1'$ satisfying $C_1'\cap C_0=\emptyset$ (hence also
$C_1'=C_1^2=e$, \cf 2.3). Then there exists $\phi\in\cal{G}(X_e,\pi)$
such that $\phi(C_1)=C_1'$. In particular, the system
$(\pi: X_e\to\PP^1;C_0,C_1)$ is uniquely determined (up to isomorphism)
by the integer $e$ and the conditions $C_0^2=-e$, $C_0\cap C_1=\emptyset$
(and from the fact that $C_0$ and $C_1$ are images of sections of
$\pi$). $C_1$ automatically satisfies $C_1^2=e$.}

\pf Let $g:\cal{O}\oplus\cal{O}(-e)\to\cal{L}$ be the surjection
corresponding to the section $\tau:\PP^1\to X_e$ with $\im\tau=C_1'$.
Then by \cite{Hartshorne}, V 2.9, $\deg\cal{L}=C_0\cdot C_1'$, hence
$\cal{L}=\cal{O}_{\PP^1}$. Now, consider the map
$g\oplus pr_2:\cal{O}\oplus\cal{O}(-e)\to\cal{O}\oplus\cal{O}(-e)$.
Since $C_0\cap C_1'=\emptyset$, this is an isomorphism, which induces
$\phi$. \edpf

\subsection{Example} (\cf \cite{Hartshorne} II 8.24.) Assume that two smooth
surfaces $E_1$ and $E_0$ in the smooth three-fold $X$ intersect each
other transversally. Additionally, assume that $E_1\cap E_0=C$ is
isomorphic to $\PP^1$ and $N_{C|E_1}=\cal{O}(a)$ and
$N_{C|E_0}=\cal{O}(b)$ with $a\le b$. Let $\widetilde{X}$ be obtained
by blowing up $X$ along $C$, and let $E$ be the exceptional divisor of
this blowing up $\pi$. Then $\pi$ induces a projection $\pi:E\to\PP^1$
making $E$ a ruled surface. In fact $E\approx X_{b-a}$.

Let $C_i$ be the intersection of $E$ with the strict transform of $E_i$
($i=0,1$). Then $C_0$ and $C_1$ are exactly the irreducible curves on
$X_{b-a}$ determined by (2.6). In particular $N_{C_1|E}=\cal{O}(b-a)$ and
$N_{C_0|E}=\cal{O}(a-b)$.

\subsection{Non-minimal ruled surfaces}

Now, we will fix a ruled surface $\pi:X_e\to\PP^1$ and a set of distinct
points $\{P_1,\dots,P_k;P_1',\dots,P_l'\}$ on $\PP^1$. Let $Q_j'$
be the unique point on $C_1$ with $\pi(Q_j')=P_j'$.

\begin{figure}
\begin{picture}(360,240)
\put(0,34){\line(1,0){340}}
\put(344,30){$\PP^1$}
\put(0,85){\line(1,0){340}}
\put(344,81){$C_1^m$}
\put(0,235){\line(1,0){340}}
\put(344,231){$C_0^m$}
\put(170,228){$-e$}
\multiput(20,31)(60,0){3}{\line(0,1){6}}
\multiput(220,31)(80,0){2}{\line(0,1){6}}
\put(18,20){$P_0$}
\put(78,20){$P_k$}
\put(45,20){$\cdots$}
\put(178,20){$\cdots$}
\put(258,20){$\cdots$}
\put(138,20){$P_{j_1}'$}
\put(120,2){{\small $m_{j_1}'=1$}}
\put(218,20){$P_{j_2}'$}
\put(198,2){{\small $m_{j_2}'=2l_{j_2}$}}
\put(298,20){$P_{j_3}'$}
\put(268,2){{\small $m_{j_3}'=2l_{j_3}+1\ge 3$}}
\qbezier(20,240)(10,200)(20,160)
\qbezier(20,160)(30,120)(20,80)
\put(24,158){$0 \ (1)$}
\qbezier(80,240)(70,200)(80,160)
\qbezier(80,160)(90,120)(80,80)
\put(84,158){$0 \ (1)$}
\qbezier(140,240)(130,200)(140,160)
\qbezier(140,160)(150,120)(140,80)
\put(144,158){$0 \ (1)$}
\qbezier(216,210)(224,225)(216,240)
\put(224,224){{\small $-1$}}
\put(240,224){{\small $(1)$}}
\qbezier(216,188)(224,203)(216,218)
\put(224,202){{\small $-2$}}
\put(240,202){{\small $(1)$}}
\qbezier(216,166)(224,181)(216,196)
\put(224,180){{\small $-2$}}
\put(240,180){{\small $(1)$}}
\qbezier(216,80)(224,95)(216,110)
\put(224,94){{\small $-1$}}
\put(240,94){{\small $(1)$}}
\qbezier(216,102)(224,117)(216,132)
\put(224,116){{\small $-2$}}
\put(240,116){{\small $(1)$}}
\qbezier(296,210)(304,225)(296,240)
\put(304,224){{\small $-1$}}
\put(320,224){{\small $(1)$}}
\qbezier(296,188)(304,203)(296,218)
\put(304,203){{\small $-2$}}
\put(320,203){{\small $(1)$}}
\qbezier(296,146)(304,161)(296,176)
\put(304,160){{\small $-2$}}
\put(320,160){{\small $(1)$}}
\qbezier(296,124)(304,139)(296,154)
\put(304,138){{\small $-3$}}
\put(320,138){{\small $(1)$}}
\qbezier(296,102)(304,117)(296,132)
\put(304,116){{\small $-1$}}
\put(320,116){{\small $(2)$}}
\qbezier(296,80)(304,95)(296,110)
\put(304,94){{\small $-2$}}
\put(320,94){{\small $(1)$}}
\put(165,74){\vector(0,-1){32}}
\put(170,55){$\pi^m$}
\end{picture}
\caption{The surface $X^m_e$ \label{sec2pg8}}
\end{figure}
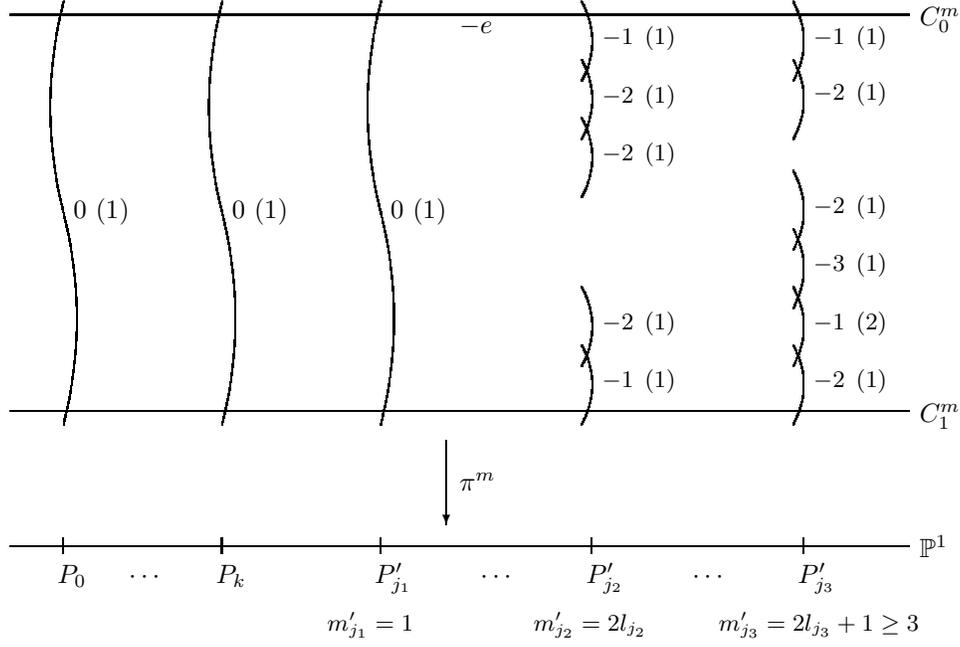

In the sequel we will modify $X_e$ by some blow ups using the following
recipe. Above the points $\{P_j\}_j$ we will not modify the fibers, we will
only mark them. On the other hand, the fibers $\pi^{-1}(P_j')$ will
(eventually) be modified. For this, we fix some integers $m_j'>0$ with
$2e=\sum_{j=1}^l m_j'$. This also means that if $e=0$, then $l=0$ and
$\{P_1',\dots,P_l'\}=\emptyset$, hence there is no modification. If
$e>0$, then we fix local coordinates $(x_j,y_j)$ in a small neighborhood
$U_j$ of $Q_j'$ with $\{x_j=0\}=C_1\cap U_j$, and we consider a local
curve $D=\{x_j^2=y_j^{m_j'}\}$ in $U_j$. Then we will modify $U_j$ by a
minimal sequence of blow ups so that the total  transform of $D\cup C_1$
form a normal crossing divisor. Here we
distinguish three cases. If $m_j'=1$, then $D=\{x_j^2=y_j\}$ intersects $C_1$
transversally, hence no blow up is needed. If $m_j'$ is odd and
$m_j'\ge 3$, then $(m_j'+3)/2$ blow ups are needed. Finally, if $m_j'$
is even, then one needs $m_j'/2$ blow ups.
We do this for any $1\le j\le l$. Notice that
the sequence of blow ups (in particular,
the analytic type of the new surface) is
independent of the choices of the local coordinates $(x_j,y_j)$ in $U_j$,
depending only on the integers $m_j'$ (and the position of the points
$\{P_j'\}_j$ in $\PP^1$).
The new surface will be denoted by
$X_e^m$.
The surface $X_e^m$ and the projection $\pi^m:X_e^m\to \PP^1$
(where $\pi^m=\pi\circ$the sequence of blow ups)
 looks as in Figure~\ref{sec2pg8}.

Here $C_i^m$ denotes the strict transform of $C_i$ ($i=0,1$).
The strict transform of $\pi^{-1}(P_j')\subset X_e$ is the unique
irreducible component of $(\pi^m)^{-1}(P_j')$ which intersects $C_0^m$.
The other components of $(\pi^m)^{-1}(P_j')$ are the new exceptional
divisors. The non-positive  integer near each irreducible curve denotes the
self-intersection of the curve in $X_e^m$.
The positive integer in parenthesis is
the multiplicity of the corresponding irreducible component in the divisor
$\pi^{\ast}(P_j')$.

It is clear that

\subsection{}
\[ \left\{
\begin{array}{rcl}
\displaystyle (C_0^m)^2&=&-e \\[12pt]
\displaystyle (C_1^m)^2&=&\displaystyle e-\sum_{m_j'\mbox{
{\footnotesize even}}}m_j'/2
   -\sum_{m_j'\mbox{ {\footnotesize odd}}\, >1}(m_j'+1)/2
\end{array}
\right.  \]

\noindent Moreover, by \cite{Hartshorne}, V 3.2, one has:

\subsection{}
$\displaystyle \Pic(X_e^m)=\Z^r, \mbox{ where }
   r=2+\sum_{m_j'\mbox{ {\footnotesize even}}}m_j'/2+
   \sum_{m_j'\mbox{{\footnotesize odd}}\, >1}(m_j'+3)/2$

\noindent and it is
generated by $C_0^m$, a generic fiber $f$, and the new irreducible
exceptional curves $\{C_j'\}$.
The intersection matrix can be easily read from
Figure~\ref{sec2pg8} using $(C_0^m)^2=-e$, $f^2=0$, $C_0^m\cdot f=1$, and
$C_0^m\cdot C_j'=f\cdot C_j'=0$ for any irreducible exceptional curve $C_j'$.

Notice that the isomorphism type of $X_e^m$ is completely determined by the
integer $e$, $\{m_j'\}_{j=1}^l$ and the position of the points $\{P_j'\}$.
On the other hand, the homeomorphism type of the
system $(X_e^m,C_0^m,C_1^m,\{(\pi^m)^{-1}(P_j)\}_j,
\{(\pi^m)^{-1}(P_j')\}_j)$ is determined completely by $e$ and $\{m_j'\}$,
and does not depend on the choice of the points $\{P_j\}_j$ and $\{P_j'\}_j$.

\subsection{}
Any  compact irreducible exceptional divisor of the birational
modification which will be constructed later, associated with
$g=f(x,y)+z^2$,  will be isomorphic to some $X_e^m$. The irreducible components
of $C_0^m\cup C_1^m\cup\left(\cup_{j}(\pi^m)^{-1}(P_j)\right)\cup
\left(\cup_{j}(\pi^m)^{-1}(P_j')\right)$
provide the intersection curves with other irreducible exceptional
divisors.

On some of the irreducible exceptional divisors $X_e^m$ we have to put one more
curve: the intersection of $X_e^m$ with the strict transform of $\{g=0\}$.
This discussion is postponed until (3.5).

\subsection{} The non-compact irreducible exceptional divisors of $g$ will be
(non-minimal) disc bundles over $\PP^1$. Notice that for any $e\in\Z$, there is
a unique disc bundle $\pi: B_e\to\PP^1$ such that the zero section $C$
satisfies $C^2=e$. Similarly as in the case of $X_e$, we can fix some data
which codifies a sequence of blow ups with centers above the
zero  section of $\pi$. In this way, we obtain the  (non-minimal)
modified disc bundle $B_e^m$ with natural projection $B_e^m\to\PP^1$. These
are the candidates for the non-compact irreducible exceptional divisors
of our birational modification.

\section{The Embedded Resolution of $g=f(x,y)+z^2$}

\subsection{} First we fix some notations regarding a convenient embedded
resolution graph of the isolated plane curve singularity
$f:(\C^2,0)\to(\C,0)$.

Let $\phi_f:Z\to U^2$ be an embedded resolution of the pair
$(f^{-1}(0),0)\subset (\C^2,0)$, where $U^2$ is a small representative of
$(\C^2,0)$. We denote the irreducible exceptional divisors by $\{A_1,\dots,
A_s\}$, the strict transforms of the irreducible components of
$\{f=0\}$ by $\{\St_1,\dots,\St_{s'}\}$, and the
collection of all these irreducible
components by $\{D_1,\dots,D_{s+s'}\}$. It is obvious
that each $A_i$ is rational. Topologically, the pair $\phi_f^{-1}(\{f=0\})
\subset Z$ (or $(\{f=0\},0)\subset(\C^2,0)$) is codified by the following data:

a) the intersection matrix $(A_i\cdot A_j)_{i,j}$;

b) the intersections $A_i\cdot\St_j$;

c) the multiplicities $m_i(f)$ of $f\circ\phi_f$ along $D_i$ (actually a, b
$\Rightarrow$ c).

\noindent In the next construction, it is convenient to assume that there is
no pair $D_i$ and $D_j$ with $D_i\cdot D_j\ne 0$ such that both $m_i(f)$ and
$m_j(f)$ are odd.

This situation can always be realized: indeed if $D_i\cap D_j\ne\emptyset$ and
$m_i(f)\cdot m_j(f)\equiv1\mod 2$, then we blow up the intersection point
$D_i\cap D_j$. The multiplicity of the new exceptional divisor will be the
even number $m_i(f)+m_j(f)$. Notice that there is a unique minimal resolution
$\phi_f$ with this property.

We will use the following notations as well. For a generic point $P$ on $A_i$,
there are local coordinates $(u,v)$ in a small neighborhood $U$ of $P$ such
that $f\circ\phi_f|U=u^{m_i(f)}$. Similarly, if $P=D_i\cap D_j$, then in some
local coordinates in a small neighborhood $U$  of $P$
one has $f\circ\phi_f|U=u^{m_i(f)} v^{m_j(f)}$.

\subsection{} In the construction of the embedded resolution graph of
$g(x,y,z)=f(x,y)+z^2$, we will use Jung's strategy. Consider the following
diagram (already mentioned above):

\[
\begin{array}{ccrcl}
\widetilde{X} &\stackrel{r}{\longrightarrow} & X &
  \stackrel{\phi'}{\longrightarrow} &
  U^3\supset g^{-1}(0) \\
 & & \Big\downarrow\vcenter{\rlap{\scriptsize $p'$}} & &
  \Big\downarrow\vcenter{\rlap{\scriptsize $p$}} \\
 & & Z & \stackrel{\phi_f}{\longrightarrow} &
  U^2 \supset f^{-1}(0)
\end{array}
\]

\noindent We have:
\begin{itemize}
\item $U^3$ is a small representative (polydisc) of $(\C^3,0)$ and $p: U^3\to
U^2$ is induced by the projection $(x,y,z)\mapsto(x,y)$.
\item $\phi_f:Z\to U^2$ is an embedded resolution of $(f^{-1}(0),0)\subset
(\C^2,0)$ as described in (3.1).
\item $p':X\to Z$ is the pull-back of $p:U^3\to U^2$ via $\phi_f$. Notice that
$X$ is smooth. Let $T_g:=(\phi')^{-1}(g^{-1}(0))$ be the total transform of
$g^{-1}(0)$ in $X$. Fix a generic point $A_i$ and small coordinate neighborhood
$U$ of $P$ as in (3.1). Then $(p')^{-1}(U)$ admits local coordinates $(u,v,z)$,
where $p'(u,v,z)=(u,v)$, and $T_g\subset(p')^{-1}(U)$ has equation $u^{m_i(f)}
+z^2=0$. Similarly, if $P=D_i\cap D_j$, then $T_g$ is given by $u^{m_i(f)}
v^{m_j(f)}+z^2=0$.
\item $r$ is an embedded resolution of $T_g\subset X$ (\cf 3.4). The composed
map $\phi'\circ r$ is denoted by $\phi_g$. By construction
$\phi_g:\widetilde{X}\to U^3$ is an ``embedded resolution'' of $(g^{-1}(0),0)
\subset(\C^3,0)$, which is an isomorphism above $U^3-\{x=y=0\}$, \cf
Introduction.
\end{itemize}

Now we will describe the  exceptional divisor $E$ of $\phi_g$.

\subsection{} The exceptional divisor $E$ is a union $E_{nc}\cup E_c$ where
$E_{nc}$ (respectively $E_c$) is the union of non-compact (respectively,
compact) irreducible exceptional components.
$E_{nc}$ is created in two steps: first we create the exceptional divisors
of $\phi'$, then we modify them by some blow ups.
Indeed, the resolution $\phi_f$ gives rise to the
exceptional curve $A=\phi_f^{-1}(0)$. This lifted, gives rise to the
exceptional surfaces $A\times D=(\phi')^{-1}(D)\subset X$, where $D$ is the
disc $\{x=y=0\}\subset U^3$. Each irreducible component of $A\times D$ has the
form $A_j\times D$, \ie it is a disc bundle with trivial self-intersection of
the zero section. We denote $A_j\times D$ by $E(A_j)$.

The multiplicity of the function $g$ (or of $z$) along each $E(A_j)$ is zero.

If $A_j'=(p'|T_g)^{-1}(A_j)\subset T_g$,
then $N_{A_j'|X}=\cal{O}\oplus\cal{O} (A_j^2)$.

\subsection{} Now we describe $r:\widetilde{X}\to X$ in more detail. First, we
fix an ordering of the irreducible exceptional divisors $\{A_j\}_{1\le j\le s}$
of $\phi_f$. By convention, if $m_i(f)$ is even and $m_j(f)$ is odd, then
$i<j$. (Sometimes we say that $A_i$ is older than $A_j$ if $i<j$.)
Such an ordering always exists and in general it is not unique.
Different orderings provide different modifications.

Now, fix an irreducible component $A_i$. Fix a generic point $P$ on $A_i$, let
$U$ be a small neighborhood of $P$ as in (3.1). Then $T_g\cap U=\{(u,v,z):
u^{m_i(f)}+z^2=0\}$, where $\{u=0\}=E(A_i)\cap U$.

The transversal plane curve singularity has type $A_{m_i(f)-1}$;
corresponding to
its minimal embedded resolution, we blow up the corresponding rational curves
above $A_i$. First we resolve this transversal singularity completely above
$A_1$, then we continue with $A_2$, and so on. If $m_i(f)$ is even, we need
$m_i(f)/2$ blow ups; if $m_i(f)=1$  then we need no modification, if
$m_i(f)$ is odd and $>1$, then we need $(m_i(f)+3)/2$ blow ups.

More precisely, assume that we finished this procedure for $A_1,\dots,A_{i-1}$,
and we want to continue with the curve $A_i$.

If $m_i(f)$ is even, $m_i(f)=2l_i$, then the graph of the
minimal embedded resolution of $(\{u^{2l_i}+z^2=0\},0)\subset (\C,0)$ is:

\bigskip

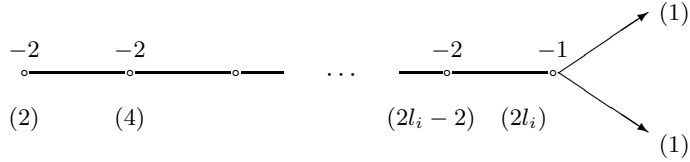
\begin{center}
\begin{picture}(290,90)(-12,0)
\multiput(0,60)(40,0){3}{\yuan}
\multiput(160,60)(40,0){2}{\yuan}
\multiput(2,60)(40,0){2}{\line(1,0){36}}
\put(162,60){\line(1,0){36}}
\put(82,60){\line(1,0){16}}
\put(114,57){$\cdots$}
\put(142,60){\line(1,0){16}}
\put(202,60){\vector(3,2){34}}
\put(202,60){\vector(3,-2){34}}
\put(-6,40){{\small $(2)$}}
\put(34,40){{\small $(4)$}}
\put(137,40){{\small $(2l_i-2)$}}
\put(180,40){{\small $(2l_i)$}}
\put(240,80){{\small $(1)$}}
\put(240,30){{\small $(1)$}}
\put(-6,66){{\small $-2$}}
\put(34,66){{\small $-2$}}
\put(154,66){{\small $-2$}}
\put(194,66){{\small $-1$}}
\end{picture}
\end{center}

\bigskip

\noindent where $(n)$ denotes the multiplicity of $u^{m_i(f)}+z^2$
along the corresponding irreducible exceptional
divisor. Corresponding to this, we make $l_i$ blow ups along the
corresponding rational curves (as axis)
above $A_i$. Let's see what types of ruled surface we will obtain.

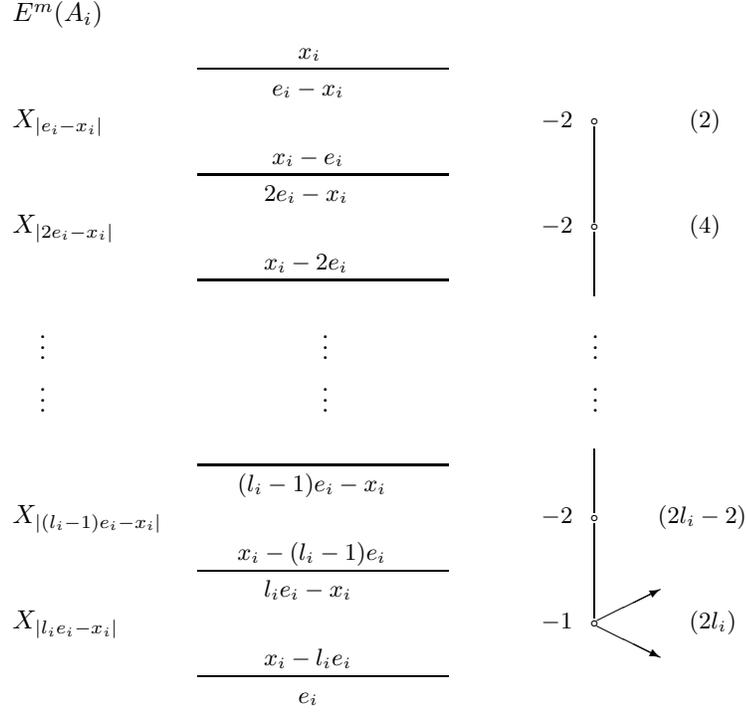
\begin{figure}
\begin{picture}(340,280)(-28,-8)
\multiput(50,0)(0,40){3}{\line(1,0){95}}
\multiput(50,150)(0,40){3}{\line(1,0){95}}
\put(88,-10){{\small $e_i$}}
\put(75,4){{\small $x_i-l_ie_i$}}
\put(75,30){{\small $l_ie_i-x_i$}}
\put(65,44){{\small $x_i-(l_i-1)e_i$}}
\put(65,70){{\small $(l_i-1)e_i-x_i$}}
\put(75,154){{\small $x_i-2e_i$}}
\put(75,180){{\small $2e_i-x_i$}}
\put(78,194){{\small $x_i-e_i$}}
\put(78,220){{\small $e_i-x_i$}}
\put(88,234){{\small $x_i$}}
\put(-20,18){$X_{|l_ie_i-x_i|}$}
\put(-20,58){$X_{|(l_i-1)e_i-x_i|}$}
\put(-20,168){$X_{|2e_i-x_i|}$}
\put(-20,208){$X_{|e_i-x_i|}$}
\put(-20,248){$E^m(A_i)$}
\multiput(200,20)(0,40){2}{\yuan}
\multiput(200,170)(0,40){2}{\yuan}
\put(200,22){\line(0,1){36}}
\put(200,62){\line(0,1){24}}
\put(200,172){\line(0,1){36}}
\put(200,168){\line(0,-1){24}}
\put(180,18){{\small $-1$}}
\put(180,58){{\small $-2$}}
\put(180,168){{\small $-2$}}
\put(180,208){{\small $-2$}}
\put(201,19){\vector(2,-1){24}}
\put(201,21){\vector(2,1){24}}
\put(236,18){{\small $(2l_i)$}}
\put(224,58){{\small $(2l_i-2)$}}
\put(236,168){{\small $(4)$}}
\put(236,208){{\small $(2)$}}
\put(-10,120){$\vdots$}
\put(97,120){$\vdots$}
\put(199,120){$\vdots$}
\put(-10,100){$\vdots$}
\put(97,100){$\vdots$}
\put(199,100){$\vdots$}
\end{picture}
\caption{Tower of ruled surfaces, the case $m_i(f)=2l_i$\label{sec3pg7}}
\end{figure}

If $A_j$ is older than $A_i$, \ie $j<i$, then $m_j(f)$ is automatically even,
hence $E(A_i)$ is modified $\sum_{j<i,A_iA_j\ne 0}m_j(f)/2$ times in different
infinitesimally close points. After these modifications
$E(A_i)$ becomes $E(A_i)'$. The strict transform of $A_i'$
in $E(A_i)'$ is denoted by $A_i''$. Then  $A_i''$ has normal bundle
$\cal{O}(-\sum_{j<i}m_j(f)/2)\oplus\cal{O}(A_i^2)$. For simplicity, we write
$A_i^2=e_i$ and $-\sum_{j<i}m_j(f)/2=x_i$.
Then we start to resolve the tranversal singularity $A_{m_i(f)-1}$ above
$A_i$.
After the sequence of blow ups along the rational curves above $A_i$,
we obtain a tower of ruled surfaces above $A_i$, as it is shown in
Figure~\ref{sec3pg7}, (\cf also with 2.7).
In this diagram the schematic picture

\begin{center}
\begin{picture}(180,60)(-26,180)
\multiput(50,190)(0,40){2}{\line(1,0){95}}
\put(94,195){{\small $n$}}
\put(88,220){{\small $-n$}}
\put(-20,208){$X_{|n|}$}
\end{picture}
\end{center}

\noindent
denotes the ruled surface $X_{|n|}$, the horizontal lines
denote the two distinguished curves $C_0$ and $C_1\subset X_{|n|}$
with self-intersections $\pm \, n$. In Figure~\ref{sec3pg7}  an
adjacency shows that $X_{|n|}$ and $X_{|m|}$ intersect each other along their
distinguished curve codified by the common horizontal line. The arrows on the
graph of $\{u^{2l_i}+z^2=0\}$ correspond to the strict transforms; their
contribution  will be discussed  in (3.5). Obviously, the integers $(n)$
denote the vanishing orders of $g$
along the corresponding components $X_{|n|}$.

Since the divisor $(f)$ of $f$ on $Z$ satisfies $(f)\cdot A_i=0$, one has
\[ e_im_i(f)+\sum_{j<i,A_jA_i\ne0}m_j(f)+\sum_{j>i,A_jA_i\ne0}m_j(f)=0. \]

Therefore, $\displaystyle
x_i-e_il_i=\frac{1}{2}\sum_{j>i,A_jA_i\ne0}m_j(f)\geq 0$.

Therefore, the surface at the bottom of the tower is $X_{x_i-e_il_i}$. The
collection $\{m_j(f)\}_{j>i,A_jA_i\ne0}$ is sometimes denoted by
$\{{m_j}'(f)\}_j$. The bottom horizontal line is the distinguished curve
$C_1$ of $X_{x_i-e_il_i}$ (since it has positive self-intersection). All the
other surfaces will be unchanged by the latter modifications, but this
$X_{x_i-e_il_i}$ will be changed by some blowing ups corresponding to the
modifications of the younger neighbors $A_j$ ($j>i$) of $A_i$.

If $m_i(f)=2l_i+1$ is odd, then the graph of the minimal resolution of
$\{u^{2l_i+1}+z^2=0\}$ is

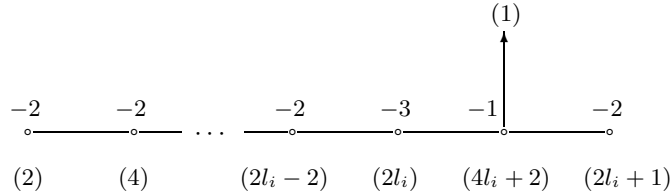
\begin{center}
\begin{picture}(270,70)(-14,34)
\multiput(0,60)(40,0){2}{\yuan}
\put(2,60){\line(1,0){36}}
\put(42,60){\line(1,0){16}}
\multiput(100,60)(40,0){4}{\yuan}
\multiput(102,60)(40,0){3}{\line(1,0){36}}
\put(180,62){\vector(0,1){36}}
\put(175,102){{\small $(1)$}}
\put(98,60){\line(-1,0){16}}
\put(63,57){$\cdots$}
\put(142,60){\line(1,0){16}}
\put(-6,40){{\small $(2)$}}
\put(34,40){{\small $(4)$}}
\put(80,40){{\small $(2l_i-2)$}}
\put(130,40){{\small $(2l_i)$}}
\put(164,40){{\small $(4l_i+2)$}}
\put(210,40){{\small $(2l_i+1)$}}
\put(-7,66){{\small $-2$}}
\put(33,66){{\small $-2$}}
\put(93,66){{\small $-2$}}
\put(133,66){{\small $-3$}}
\put(166,66){{\small $-1$}}
\put(213,66){{\small $-2$}}
\end{picture}
\end{center}

\begin{figure}
\begin{picture}(360,360)(-28,-98)
\multiput(50,-80)(0,40){5}{\line(1,0){95}}
\multiput(50,150)(0,40){3}{\line(1,0){95}}
\put(88,-90){{\small $e_i$}}
\put(65,-76){{\small $x_i-(l_i+1)e_i$}}
\put(65,-50){{\small $(l_i+1)e_i-x_i$}}
\put(60,-36){{\small $2x_i-(2l_i+1)e_i$}}
\put(60,-10){{\small $(2l_i+1)e_i-2x_i$}}
\put(75,4){{\small $x_i-l_ie_i$}}
\put(75,30){{\small $l_ie_i-x_i$}}
\put(65,44){{\small $x_i-(l_i-1)e_i$}}
\put(65,70){{\small $(l_i-1)e_i-x_i$}}
\put(78,194){{\small $x_i-e_i$}}
\put(78,220){{\small $e_i-x_i$}}
\put(88,234){{\small $x_i$}}
\put(-20,-62){$X_{|(l_i+1)e_i-x_i|}$}
\put(-20,-22){$X_{|(2l_i+1)e_i-2x_i|}$}
\put(-20,18){$X_{|l_ie_i-x_i|}$}
\put(-20,208){$X_{|e_i-x_i|}$}
\put(-20,248){$E^m(A_i)$}
\multiput(200,-60)(0,40){4}{\yuan}
\multiput(200,170)(0,40){2}{\yuan}
\multiput(200,-58)(0,40){3}{\line(0,1){36}}
\put(200,62){\line(0,1){24}}
\put(200,172){\line(0,1){36}}
\put(200,168){\line(0,-1){24}}
\put(180,-62){{\small $-2$}}
\put(180,-22){{\small $-1$}}
\put(180,18){{\small $-3$}}
\put(180,58){{\small $-2$}}
\put(180,168){{\small $-2$}}
\put(180,208){{\small $-2$}}
\put(202,-20){\vector(1,0){28}}
\put(224,-62){{\small $(2l_i+1)$}}
\put(236,-22){{\small $(4l_i+2)$}}
\put(236,18){{\small $(2l_i)$}}
\put(224,58){{\small $(2l_i-2)$}}
\put(236,168){{\small $(4)$}}
\put(236,208){{\small $(2)$}}
\put(-10,120){$\vdots$}
\put(97,120){$\vdots$}
\put(199,120){$\vdots$}
\put(-10,100){$\vdots$}
\put(97,100){$\vdots$}
\put(199,100){$\vdots$}
\end{picture}
\caption{Tower of ruled surfaces, the case $m_i(f)=2l_i+1$\label{sec3pg10}}
\end{figure}
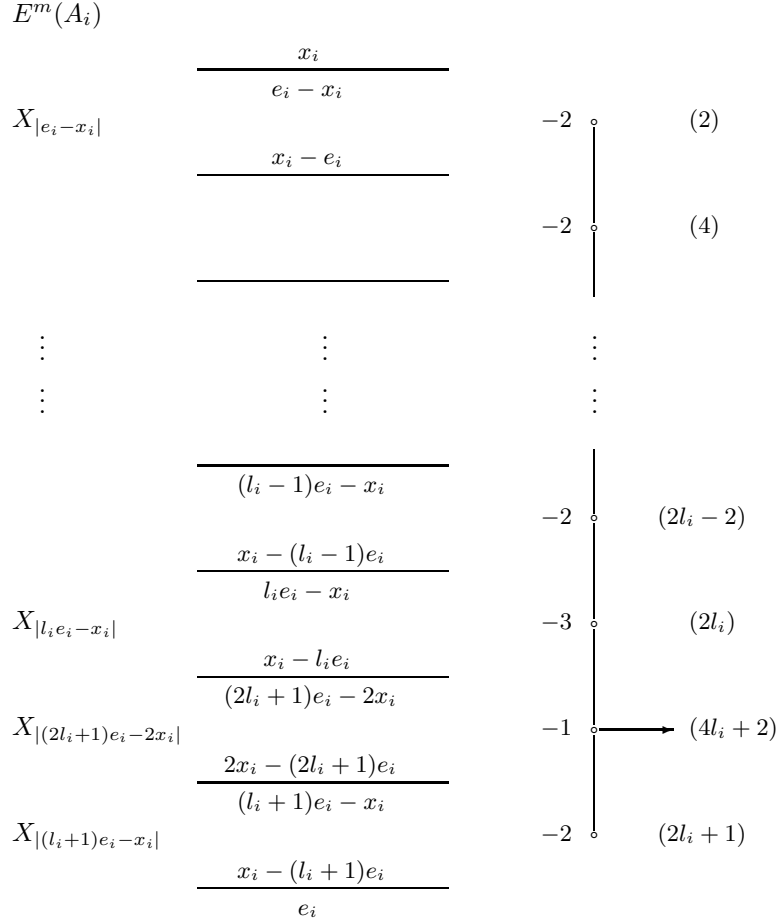

\noindent
Since all the neighbors of $A_i$ are older than $A_i$ (since all of them have
even multiplicity), $N_{A_i'|X}=\cal{O}(x_i)\oplus\cal{O}(e_i)$ where
$x_i=-\sum_{i\ne j,A_iA_j\ne 0}m_j(f)/2$ and $e_i=A_i^2$ as above. Then the
tower of ruled surfaces is as in Figure~\ref{sec3pg10}.

The relation $(f)\cdot A_i=0$ now reads as
\[ m_i(f)e_i+\sum_{A_jA_i\ne 0,i\ne j}m_j(f)=0 \]
\ie $(2l_i+1)e_i-2x_i=0$. Therefore, the second surface from the bottom (which
will support the strict transform corresponding to the arrow) is
$X_0=\PP^1\times\PP^1$.

Now we will analyze how we have to glue all these towers
of surfaces. This discussion will clarify also how we
have to modify $X_{x_i-e_il_i}$ (the case $m_i(f)$ even)
corresponding to the blowing ups of the younger neighbors.

\begin{figure}
\begin{picture}(360,360)
\put(40,0){\line(1,0){280}}
\multiput(40,360)(5,0){56}{\line(1,0){3}}
\put(160,0){\line(0,1){360}}
\put(159,330){\line(0,1){30}}
\put(164,342){$\widetilde{D}$}
\multiput(160,30)(0,30){2}{\line(1,0){150}}
\multiput(160,105)(0,30){3}{\line(1,0){150}}
\multiput(40,195)(0,30){2}{\line(1,0){120}}
\multiput(40,270)(0,30){3}{\line(1,0){120}}
\put(144,12){{\small $-1$}}
\put(144,42){{\small $-2$}}
\put(144,117){{\small $-2$}}
\put(144,147){{\small $-2$}}
\put(144,177){{\small $-1$}}
\put(150,207){{\small $0$}}
\put(150,282){{\small $0$}}
\put(150,312){{\small $0$}}
\multiput(164,12)(0,30){2}{{\small $0$}}
\multiput(164,117)(0,30){2}{{\small $0$}}
\put(164,177){{\small $-1$}}
\put(164,207){{\small $-2$}}
\multiput(164,282)(0,30){2}{{\small $-2$}}
\put(80,170){{\small $(2t_j)$}}
\put(70,207){{\small $(2t_j-2)$}}
\put(84,282){{\small $(4)$}}
\put(84,312){{\small $(2)$}}
\put(84,342){{\small $(0)$}}
\put(6,90){$X^m_{x_j-t_je_j}$}
\put(6,312){$X_{|e_j-x_j|}$}
\put(6,342){$E^m(A_j)$}
\put(88,242){$\vdots$}
\put(234,12){{\small $(2l_i)$}}
\put(227,42){{\small $(2l_i-2)$}}
\put(236,117){{\small $(4)$}}
\put(236,147){{\small $(2)$}}
\put(236,227){{\small $(0)$}}
\put(227,312){$E^m(A_i)$}
\put(320,12){$X^m_{x_i-l_ie_i}$}
\put(320,147){$X^m_{|e_i-x_i|}$}
\put(240,78){$\vdots$}
\end{picture}
\caption{Gluing, the case $m_i(f)=2l_i$, $m_j(f)=2t_j$ and
$j<i$\label{sec3pg12}}
\end{figure}

If $P$ is an intersection point $A_i\cap A_j$, then $T_g\subset(p')^{-1}(U)$
has local equation $u^{m_i(f)}v^{m_j(f)}+z^2=0$. We distinguish two cases.
In the first case $m_i(f)=2l_i$,
$m_j(f)=2t_j$, and we ssume that the $u$-axis is older than the $v$-axis
(\ie $j<i$). Then
above $U$, the exceptional divisor $E$ has the form shown as in
Figure~\ref{sec3pg12}. In the picture, the $\widetilde{D}$ is a disc over
$D=\{x=y=0\}$. The vertical segments of
contact between projective surfaces denote rational curves. The integers
$a|b$ denote the self-intersection numbers of this curve in the two surfaces
correspondingly.

\begin{figure}
\begin{picture}(360,390)(0,-30)
\put(40,-30){\line(1,0){280}}
\multiput(40,360)(5,0){56}{\line(1,0){3}}
\put(160,-30){\line(0,1){390}}
\put(159,330){\line(0,1){30}}
\put(164,342){$\widetilde{D}$}
\multiput(160,0)(0,30){3}{\line(1,0){150}}
\multiput(160,105)(0,30){3}{\line(1,0){150}}
\multiput(40,195)(0,30){2}{\line(1,0){120}}
\multiput(40,270)(0,30){3}{\line(1,0){120}}
\put(144,-18){{\small $-2$}}
\put(144,12){{\small $-1$}}
\put(144,42){{\small $-3$}}
\put(144,117){{\small $-2$}}
\put(144,147){{\small $-2$}}
\put(144,177){{\small $-1$}}
\put(150,207){{\small $0$}}
\put(150,282){{\small $0$}}
\put(150,312){{\small $0$}}
\multiput(164,-18)(0,30){3}{{\small $0$}}
\multiput(164,117)(0,30){2}{{\small $0$}}
\put(164,177){{\small $-1$}}
\put(164,207){{\small $-2$}}
\multiput(164,282)(0,30){2}{{\small $-2$}}
\put(80,170){{\small $(2t_j)$}}
\put(70,207){{\small $(2t_j-2)$}}
\put(84,282){{\small $(4)$}}
\put(84,312){{\small $(2)$}}
\put(84,342){{\small $(0)$}}
\put(6,90){$X^m_{x_j-t_je_j}$}
\put(6,312){$X_{|e_j-x_j|}$}
\put(6,342){$E^m(A_j)$}
\put(88,242){$\vdots$}
\put(227,-18){{\small $(2l_i+1)$}}
\put(227,12){{\small $(4l_i+2)$}}
\put(234,42){{\small $(2l_i)$}}
\put(236,117){{\small $(4)$}}
\put(236,147){{\small $(2)$}}
\put(236,227){{\small $(0)$}}
\put(227,312){$E^m(A_i)$}
\put(320,147){$X^m_{|e_i-x_i|}$}
\put(240,78){$\vdots$}
\end{picture}
\caption{Gluing, the case $m_i(f)=2l_i+1$, $m_j(f)=2t_j$ \label{sec3pg13}}
\end{figure}

If the local equation is $u^{m_i(f)}v^{m_j(f)}+z^2=0$, where $m_i(f)=2l_i+1$,
$m_j(f)=2t_j$, then automatically $j<i$, and
$E$ above a neighborhood of $P$ looks as in Figure~\ref{sec3pg13}.

Now, it is possible to verify using a  local equation  of type
$u^{m_i(f)}v^{m_j(f)}+z^2=0$, that in both cases after the steps
described above the total transform  of $\{g=0\}$ is a normal crossing
divisor, \ie we do not have to blow up any other center, and the resolution
procedure ends (\cf Orbanz \cite{Orbanz}.)

This ends the complete description of the exceptional set $E$ and of all the
normal bundles $N_{E_k\cap E_{k'}|E_k}$ ($k\ne k'$).

\subsection{The intersection of the strict transform $\St(g)$ of $\{g=0\}$
with $E$}
Corresponding to the arrows with multiplicity $(1)$
of the embedded resolution graph of the plane curve
singularity $\{u^m+z^2=0\}$, in the tower of surfaces
above any $A_i$, there is exactly one, say $\widetilde{E}(A_i)$, which
intersects the strict transform $\St(g)$. Let us denote this intersection by
$S_i^m=\St(g)\cap\widetilde{E}(A_i)$.

First we will describe the position of $S_i^m$ in $\widetilde{E}(A_i)$.

If $m_i(f)$ is odd, the $\widetilde{E}(A_i)=X_{(2l_i+1)e_i-2x_i}=X_0=
\PP^1\times\PP^1$. If by this identification $\pi:\widetilde{E}(A_i)\to A_i$
corresponds to the first projection, and $C_0$ and $C_1$, the intersection
curves with the neighbor
ruled surfaces in the same tower correspond to two fibers of
the second projection, then $S_i^m$ provides another section of $\pi$ with
$S_i^m\cap C_0=S_i^m\cap C_1=\emptyset$. Therefore
 $S_i^m$ is another fiber of the
second projection. Obviously its self-intersection is $(S_i^m)^2=0$ and
$S_i^m\cdot f=1$.

The situation is slightly more complicated if $m_i(f)=2l_i$ is even.

Consider that moment of the resolution procedure when we finished the
construction of the tower about $A_i$:
we just created  the last ruled surface
$X_{x_i-l_ie_i}$, but we did not start the next tower above $A_{i+1}$.
Consider the intersection points $\{P_0,\dots,P_k\}$ of $A_i$ with older
exceptional curves $A_j$ ($j<i$), and also all the
intersection points $\{P_1',\dots, P_l'\}$ of $A_i$ with all the other
irreducible components of the total transform of $\{f=0\}$, i.e. with the
younger exceptional curves $A_j$ ($j>i$) and with the
strict transforms $\St_j$.
Then, similarly as above , we denote the collection of multiplicities
$m_j(f)$ of $f$ along the components $A_j \ (A_j\cdot A_i\not=0, \ j>i)$
and $\St_j \ (\St_j\cdot A_i\not=0)$ by $\{m_t'\}^l_{t=1}$ (such that the
index corresponds to the index of the points $P_t'$).  Then we are in the
situation of (2.8) where $X_e=X_{x_i-l_ie_i}$.
Using the notation of (2.8), the modified surface $X_e^m$ is exactly the
surface which is obtained from  $X_{x_i-l_ie_i}$ after we finish the
resolution procedure, \ie we construct all  the neighbor towers as well.

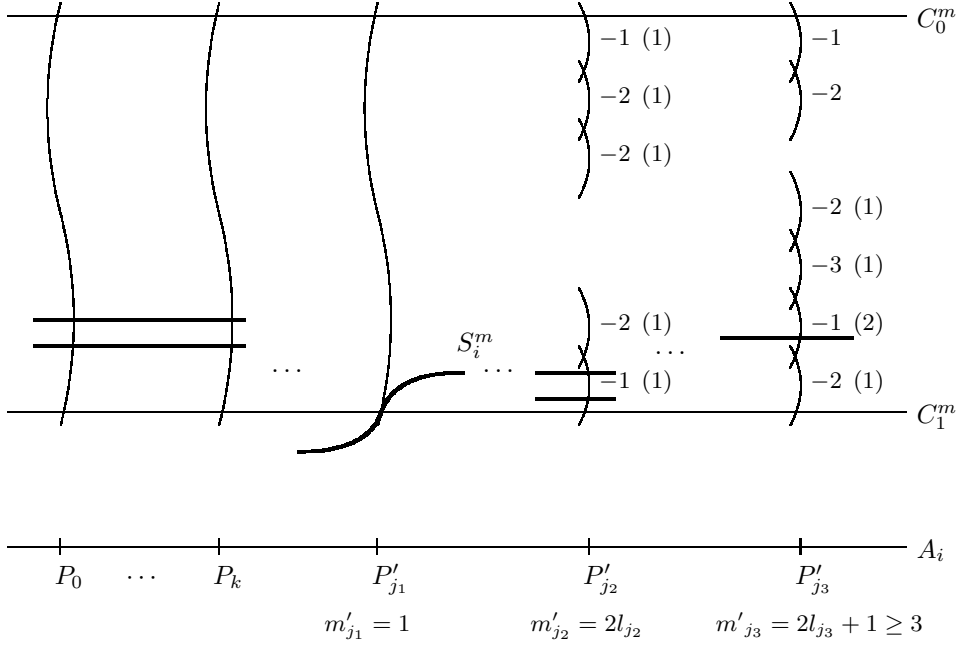
\begin{figure}
\begin{picture}(360,240)
\put(0,34){\line(1,0){340}}
\put(344,30){$A_i$}
\put(0,85){\line(1,0){340}}
\put(344,81){$C_1^m$}
\put(0,235){\line(1,0){340}}
\put(344,231){$C_0^m$}
\multiput(20,31)(60,0){3}{\line(0,1){6}}
\multiput(220,31)(80,0){2}{\line(0,1){6}}
\put(18,20){$P_0$}
\put(78,20){$P_k$}
\put(45,20){$\cdots$}
\put(138,20){$P_{j_1}'$}
\put(120,2){{\small $m'_{j_1}=1$}}
\put(218,20){$P'_{j_2}$}
\put(198,2){{\small $m'_{j_2}=2l_{j_2}$}}
\put(298,20){$P'_{j_3}$}
\put(268,2){{\small ${m'}_{j_3}=2l_{j_3}+1\ge 3$}}
\qbezier(20,240)(10,200)(20,160)
\qbezier(20,160)(30,120)(20,80)
\qbezier(80,240)(70,200)(80,160)
\qbezier(80,160)(90,120)(80,80)
\qbezier(140,240)(130,200)(140,160)
\qbezier(140,160)(150,120)(140,80)
\qbezier(216,210)(224,225)(216,240)
\put(224,224){{\small $-1$}}
\put(240,224){{\small $(1)$}}
\qbezier(216,188)(224,203)(216,218)
\put(224,202){{\small $-2$}}
\put(240,202){{\small $(1)$}}
\qbezier(216,166)(224,181)(216,196)
\put(224,180){{\small $-2$}}
\put(240,180){{\small $(1)$}}
\qbezier(216,80)(224,95)(216,110)
\put(224,94){{\small $-1$}}
\put(240,94){{\small $(1)$}}
\qbezier(216,102)(224,117)(216,132)
\put(224,116){{\small $-2$}}
\put(240,116){{\small $(1)$}}
\qbezier(296,210)(304,225)(296,240)
\put(304,224){{\small $-1$}}
\qbezier(296,188)(304,203)(296,218)
\put(304,203){{\small $-2$}}
\qbezier(296,146)(304,161)(296,176)
\put(304,160){{\small $-2$}}
\put(320,160){{\small $(1)$}}
\qbezier(296,124)(304,139)(296,154)
\put(304,138){{\small $-3$}}
\put(320,138){{\small $(1)$}}
\qbezier(296,102)(304,117)(296,132)
\put(304,116){{\small $-1$}}
\put(320,116){{\small $(2)$}}
\qbezier(296,80)(304,95)(296,110)
\put(304,94){{\small $-2$}}
\put(320,94){{\small $(1)$}}
\linethickness{1pt}
\qbezier(110,70)(136,70)(141,85)
\qbezier(172,100)(146,100)(141,85)

\put(200,100){\line(1,0){30}}
\put(200,90){\line(1,0){30}}
\put(10,120){\line(1,0){80}}
\put(270,113){\line(1,0){50}}
\put(10,110){\line(1,0){80}}
\put(100,100){$\ldots$}
\put(180,100){$\ldots$}
\put(245,107){$\ldots$}
\put(170,108){$S_i^m$}
\end{picture}
\caption{$S^m_i$, the case $m_i(f)=2l_i$\label{sec3pg17}}
\end{figure}

Moreover, the intersection with the strict transform can be identified
as follows.
Assume that the $x$-chart $\C_x$ of $A_i$ contains all the points
$\{P_t'\}^l_{t=1}$. Then the intersection $S_i$ of the strict transform
of $g$  with $X_e$ in the $\C_x\times\PP^1$ chart is
\[ \left\{(x,[u_0:u_1]): u_0^2=u_1^2\prod_{t=1}^l(x-x_t)^{m_t'}\right\}. \]
In particular, it is uniquely determined by the pair
$(A_i, \{P_t'\}_t)$
and from the numerical data. After we blow up $X_e$ and we obtain $X_e^m$,
denote the
strict transform of $S_i$ in $X_e^m$ by
$S_i^m=\St(g)\cap X_e^m$. Schematically, $S_i^m$
is as in Figure~\ref{sec3pg17} (but this is a ``real picture'', which does not
reflect exactly the ``complex picture''). In the picture, the points $P_i'$
with $m_i'=1$  correspond exactly to the intersection points  of $A_i$
with the components of the strict transform   of $\{f=0\}$.
Above these points $S_i^m$ intersects transversally $C_1^m$, and these are the
only intersection points of $S_i^m$ and  $C_1^m$.

Therefore,
\begin{eqnarray*}
 S_i^m\cdot C_1^m&=&
     \#\{\text{strict transforms of $f=0$ supported on $A_i$}\} \\
 &=&\#\{m_j(f)=1; D_j\cdot A_i\ne 0\}
\end{eqnarray*}

Since $S_i^m\to A_i$ is a double covering with branch locus exactly over the
points $P_j'$ with $m_j'$ odd, one has:
\begin{itemize}
\item $S_i^m\cdot f=2$
\item If each ${m_j}'$ is even (\ie $A_i$ has no neighbor with odd
multiplicity) then $S_i^m$ has two disjoint irreducible components, both
isomorphic to $A_i(\approx\PP^1)$.
\item If
at least one $m_j'$ is odd, then $S_i^m$ is irreducible; its genus can be
computed by an Euler characteristic argument (or Hurwitz's formula, see
\cite{Hartshorne} page 299):
\[ \text{genus}(S_i^m)=\frac{\#\{j: m_j' \text{ odd}\}-2}{2}. \]
\end{itemize}
Now we determine the self-intersection of $S_i^m$.
The above diagram shows that $S_i^m\cdot C_0^m=0$ and also  one can read
all the intersections
of $S_i^m$ with the new exceptional divisors of $X_e^m$.

If  $m_j'$ is  odd and $\ge 3$, then above $P_j'$
(i.e. in $(\pi^m)^{-1}(P_j')$)
 we distinguish two irreducible
curves $F_j'$ and $F_j$  defined by $F_j\cdot C_i^m=1$, $F_j'\cdot F_j=1$.

Notice that $F_j'$ is the unique curve with multiplicity two in the divisor
$(\pi^m)^{\ast}(P_j')$ and $S_i^m$ intersects the fiber $(\pi^m)^{-1}(P_j')$
exactly along $F_j'$. Therefore $S_i^m\cdot F_j'=1$. We invite the reader
to verify that in $\Pic(X_e^m)$ one has
\[ S_i^m=2C_1^m+\sum F_j, \quad (\text{sum over $m_j'$ odd $\ge3$}). \]
Therefore,
\begin{eqnarray*}
(S_i^m)^2&=&(2C_1^m+\sum F_j)^2 \\
 &=&4(C_i^m)^2+4M_i-2M_i=4(C_i^m)^2+2M_i,
\end{eqnarray*}
where $M_i=\#\{j:m_j' \text{ odd},m_j'\ge 3\}$.

If $m_i(f)$ is even, and there is no $D_i$ with $D_i\cdot A_i\ne 0$ with odd
multiplicity, then $S_i^m$ has two disjoint components, each $\approx C_1^m
\approx \PP^1$ and $\equiv C_1^m$ in $\Pic(X_e^m)$.

\subsection{The self-intersections $E_{k}^2$.}
Consider all the compact irreducible exceptional divisors  
$\{E_k\}_k$ of the resolution $\widetilde{\phi}$. In this 
paragraph we determine the ``self-intersections'' $E_k^2:=
{\cal O}_{\widetilde{X}}(E_k)|E_k\in Pic(E_k)$. 
First notice that in  the
previous discussions we have determined completely the divisor 
$(g\circ\widetilde{\phi})=St(g)+\sum_l\ m_l E_l$
(where the sum is over the compact  irreducible exceptional divisors).
But, for any $k$, 
$(g\circ\widetilde{\phi})\cdot E_k=0$ in  $Pic(E_k)$, therefore:
$$m_k E_k^2=(-St(g)-\sum_{l \not=k }\ m_l E_l )E_k$$
in $Pic(E_k)$. But all the intersections $St(g) E_k$
and $E_l E_k$ ($l\not=k$) are determined above, hence $E_k^2$ follows.

\subsection{The resolution of $(\{f(x,y)+z^2=0\},0)$}
Notice that $\widetilde{\phi}:\widetilde{X}\to\C^3$
induces a map $\widetilde{\phi}:\St(g)\to
\{g=0\}$, which is a resolution of the normal surface singularity $\{g=0\}$.
Its exceptional curve is exactly $\bigcup_iS_i^m$. In this subsection we
determine the self-intersections of the irreducible components of
$\bigcup_iS_i^m$, and their combinatorics. In particular, we re-obtain the dual
resolution graph of this resolution (which was known \cf \cite{Laufer},
see also \cite{Nemethi}). The details are left to the reader.

We distinguish several cases.

If $m_i(f)$ is odd, then $\pi:S_i^m\to A_i$ is an isomorphism, hence $S_i^m$
is rational. Above each intersection point $A_i\cap A_j$ we have exactly one
intersection point $S_i^m\cap S_j^m$. The self-intersection of $S_i^m$ in
$\St(g)$ is $e_i/2$ where $e_i=A_i^2$.

If $m_i(f)$ is even, and any $D_j$ with $D_j\cdot A_i\ne 0$ has even
multiplicity, then $S_i^m$ has two irreducible components, each isomorphic
to $A_i=\PP^1$. In this case the self-intersection of each component in
$\St(g)$ is $e_i$. Above an intersection point $A_i\cap A_j$ we have exactly
two intersection points of $S_i^m\cap S_j^m$.

If $m_i(f)$ is even, and at least one
 of the neighbors (including the strict transform components)
has  odd multiplicity $m_j(f)$, then $S_i^m$ is an irreducible curve with
$\text{genus}=(\lambda-2)/2$ where $\lambda$ is the number of odd neighbors
(including the strict transform components) (\cf also with
3.5). Its self-intersection
is $2e_i$. Above an intersection point $A_i\cap D_j$ with $m_j(f)$ even, there
are two intersection points; if $m_j(f)$ is odd, then only one intersection
point of $S_i^m\cap S_j^m$.

This determines completely the dual resolution graph of $\{f(x,y)+z^2=0\}$
from the embedded resolution graph of $f$. If we start with the minimal
embedded resolution of $f$ which has the additional property that
there are no neighbors $D_i$ and $D_j$, both with odd multiplicity,
then the constructed resolution $St(g)\to \{g=0\}$ is exactly the
{\em canonical} resolution of $\{g=0\}$ (\cf \cite{Laufer}).

The above statements about the self-intersections of the components of $S_i^m$
in $\St(g)$ can be obtained by the following ``triple point formula'' as
well.

Let $h:(\C^3,0)\to(\C,0)$ be the germ of an analytic function and let $(h)$ be
the divisor of $h\circ\widetilde{\phi}$.
We assume that $(h)$ is a normal crossing divisor.
We write $(h)=\sum_Dm_DD$ where the sum runs over the irreducible exceptional
divisors of $\widetilde{\phi}$
and the irreducible components of the strict transform of
$h$. Let $C$ be a compact curve determined by the intersection of two
components $D_1$ and $D_2$. Then
\[ m_{D_1}(C^2 \text{ in } D_2)+m_{D_2}(C^2 \text{ in } D_2)+\sum m_D=0 \]
where the last sum is over the triple points $D\cap D_1\cap D_2$.
(For a proof, see e.g. \cite{Persson}.)

In our case, in order to obtain the self-intersections of $S_i^m$'s in
$\St(g)$, we apply the above relation for $h=g=f+z^2$.

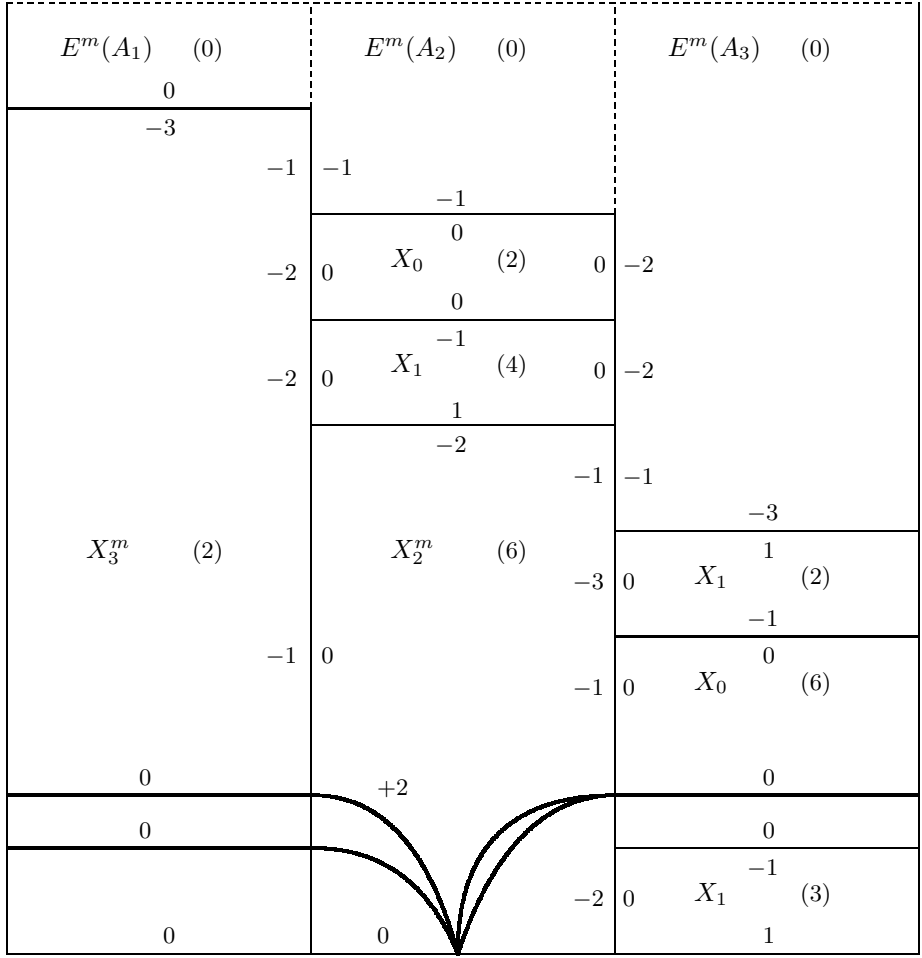
\begin{figure}
\begin{picture}(345,360)
\multiput(0,0)(345,0){2}{\line(0,1){360}}
\put(115,0){\line(0,1){320}}
\put(230,0){\line(0,1){280}}
\put(98,110){{\small $-1$}}
\put(98,215){{\small $-2$}}
\put(98,255){{\small $-2$}}
\put(98,295){{\small $-1$}}
\put(119,110){{\small $0$}}
\put(119,215){{\small $0$}}
\put(119,255){{\small $0$}}
\put(119,295){{\small $-1$}}
\multiput(115,320)(0,4){10}{\line(0,1){2}}
\multiput(230,280)(0,4){20}{\line(0,1){2}}
\put(214,18){{\small $-2$}}
\put(214,98){{\small $-1$}}
\put(214,138){{\small $-3$}}
\put(214,178){{\small $-1$}}
\put(222,218){{\small $0$}}
\put(222,258){{\small $0$}}
\put(233,18){{\small $0$}}
\put(233,98){{\small $0$}}
\put(233,138){{\small $0$}}
\put(233,178){{\small $-1$}}
\put(233,218){{\small $-2$}}
\put(233,258){{\small $-2$}}
\put(162,190){{\small $-2$}}
\put(168,203){{\small $1$}}
\put(162,230){{\small $-1$}}
\put(168,244){{\small $0$}}
\put(168,270){{\small $0$}}
\put(162,283){{\small $-1$}}
\put(286,4){{\small $1$}}
\put(280,30){{\small $-1$}}
\put(286,44){{\small $0$}}
\put(286,110){{\small $0$}}
\put(280,124){{\small $-1$}}
\put(286,150){{\small $1$}}
\put(280,164){{\small $-3$}}
\put(0,0){\line(1,0){345}}
\multiput(0,360)(4,0){86}{\line(1,0){2}}
\put(0,320){\line(1,0){115}}
\multiput(115,200)(0,40){3}{\line(1,0){115}}
\put(230,40){\line(1,0){115}}
\multiput(230,120)(0,40){2}{\line(1,0){115}}
\put(59,4){{\small $0$}}
\put(140,4){{\small $0$}}
\put(140,60){{\small $+2$}}
\put(52,310){{\small $-3$}}
\put(59,324){{\small $0$}}
\put(30,150){$X^m_3$}
\put(70,150){{\small $(2)$}}
\put(20,340){$E^m(A_1)$}
\put(70,340){{\small $(0)$}}
\put(135,340){$E^m(A_2)$}
\put(185,340){{\small $(0)$}}
\put(250,340){$E^m(A_3)$}
\put(300,340){{\small $(0)$}}
\put(145,260){$X_0$}
\put(185,260){{\small $(2)$}}
\put(145,220){$X_1$}
\put(185,220){{\small $(4)$}}
\put(145,150){$X^m_2$}
\put(185,150){{\small $(6)$}}
\put(260,20){$X_1$}
\put(300,20){{\small $(3)$}}
\put(260,100){$X_0$}
\put(300,100){{\small $(6)$}}
\put(260,140){$X_1$}
\put(300,140){{\small $(2)$}}
\linethickness{1pt}
\qbezier(170,0)(154,40)(115,40)
\put(115,40){\line(-1,0){115}}
\put(50,44){{\small $0$}}
\qbezier(170,0)(154,60)(115,60)
\put(115,60){\line(-1,0){115}}
\put(50,64){{\small $0$}}
\qbezier(170,0)(170,60)(230,60)
\put(230,60){\line(1,0){115}}
\put(286,64){{\small $0$}}
\qbezier(170,0)(190,60)(230,60)
\end{picture}
\caption{The exceptional divisor, case $x^2+y^3+z^2$\label{sec4pg1}}
\end{figure}

\subsection{Example} Assume that $f(x,y)=x^2+y^3$.
We start with the following resolution of $f$.

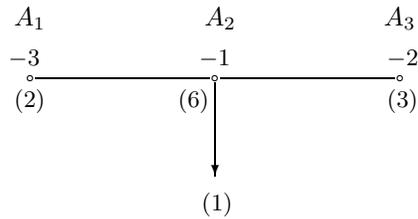
\begin{center}
\begin{picture}(160,92)
\multiput(10,60)(70,0){3}{\yuan}
\multiput(12,60)(70,0){2}{\line(1,0){66}}
\put(80,58){\vector(0,-1){35}}
\put(75,10){{\small $(1)$}}
\put(4,49){{\small $(2)$}}
\put(66,49){{\small $(6)$}}
\put(145,49){{\small $(3)$}}
\put(2,65){{\small $-3$}}
\put(74,65){{\small $-1$}}
\put(145,65){{\small $-2$}}
\put(4,80){$A_1$}
\put(76,80){$A_2$}
\put(144,80){$A_3$}
\end{picture}
\end{center}

\noindent We fix the order $A_1$,
$A_2$, $A_3$ as it is indicated above. Then $x_1=0$, $e_1=-3$; $x_2=-1$,
$e_2=-1$; $x_3=-3$, $e_3=-2$. For the components of $E$,
See Figure~\ref{sec4pg1}.



The dual resolution graph of $\displaystyle\bigcup_iS_i^m\subset\St(g)$ is

\begin{center}
\begin{picture}(160,116)
\multiput(80,60)(70,0){2}{\yuan}
\put(82,60){\line(1,0){66}}
\put(16,28){\yuan}
\put(16,92){\yuan}
\put(17,29){\line(2,1){61}}
\put(17,91){\line(2,-1){61}}
\put(0,24){{\small $-3$}}
\put(0,88){{\small $-3$}}
\put(74,65){{\small $-2$}}
\put(145,65){{\small $-1$}}
\put(8,0){$S^m_1$}
\put(76,0){$S^m_2$}
\put(144,0){$S^m_3$}
\end{picture}
\end{center}


\nocite{*}
\bibliographystyle{amsplain}
\bibliography{amsl-bib}

\end{document}